\documentclass[12pt]{article}
\usepackage[latin1]{inputenc}
\usepackage{amsmath}
\usepackage{amssymb,amsfonts}
\usepackage{graphicx}
\usepackage{times,amssymb,amscd}
\usepackage[USenglish]{babel}
\usepackage{verbatim}
\usepackage{enumerate}

\newcommand{\bC}{\mathbf{C}}

\newcommand{\bG}{\mathbf{G}}

\newcommand{\bb}{\mathbf{b}}
\newcommand{\bw}{\mathbf{w}}
\newcommand{\bv}{\mathbf{v}}

\newcommand{\bg}{\mathbf{g}}
\newcommand{\bT}{\mathbf{T}}

\newcommand{\bu}{\mathbf{u}}

\newcommand{\Bu}{\boldsymbol{u}}

\newcommand{\cD}{\mathcal{D}}

\newcommand{\cC}{\mathcal{C}}
\newcommand{\cB}{\mathcal{B}}

\usepackage[usenames,dvipsnames]{xcolor} 
\usepackage{tikz, tikz-3dplot, pgfplots}
\usepackage{tkz-graph}
\usetikzlibrary[positioning,patterns] 

\numberwithin{equation}{section}

\usepackage{algorithm}

\begin{document}
\pagestyle{myheadings}
\markboth{\centerline{R.~T.~Kozma and J.~Szirmai}}
{Structure and Visualization of Optimal Horoball Packings}
\title
{Structure and Visualization of Optimal Horoball Packings in $3$-dimensional Hyperbolic Space
}

\author{\medbreak \medbreak {\normalsize{}} \\
\normalsize Robert T. Kozma$^{(1),(2)}$ and Jen\H{o}  Szirmai$^{(1)}$ \\
\normalsize (1) Budapest University of Technology and Economics\\
\normalsize Institute of Mathematics, Department of Geometry \\
\normalsize H-1521 Budapest, Hungary \\
\normalsize (2) Department of Mathematics, Statistics, and Computer Science\\
\normalsize University of Illinois at Chicago \\
\normalsize Chicago IL 60607 USA \\
\normalsize E-mail: rkozma2@uic.edu,~szirmai@math.bme.hu }

\date{\normalsize (\today)}

\maketitle

\begin{abstract}                                       %

Four packings of hyperbolic 3-space are known to yield the optimal packing density of $0.85328\dots$. 
They are realized in the regular tetrahedral and cubic Coxeter honeycombs with Schl\"afli symbols $\{3,3,6 \}$ and $\{4,3,6\}$.  
These honeycombs are totally asymptotic, and the packings consist of horoballs (of different types) centered at the ideal vertices. 
We describe a method to visualize regular horoball packings 
of extended hyperbolic 3-space $\overline{\mathbb{H}}^3$ using the Beltrami-Klein model and the Coxeter group of the packing. We produce the first known images of these four optimal horoball packings. 

\end{abstract}

\newtheorem{theorem}{Theorem}[section]
\newtheorem{corollary}[theorem]{Corollary}
\newtheorem{lemma}[theorem]{Lemma}
\newtheorem{exmple}[theorem]{Example}
\newtheorem{defn}[theorem]{Definition}
\newtheorem{rmrk}[theorem]{Remark}
\newtheorem{proposition}[theorem]{Proposition}
\newenvironment{definition}{\begin{defn}\normalfont}{\end{defn}}
\newenvironment{remark}{\begin{rmrk}\normalfont}{\end{rmrk}}
\newenvironment{example}{\begin{exmple}\normalfont}{\end{exmple}}
\newenvironment{acknowledgement}{Acknowledgement}


\section{Introduction}

Images of the hyperbolic plane are abundant, one need not search long to find a myriad of hyperbolic plane kaleidoscopes. 
Most readers are no doubt are familiar with the work of M. C. Escher whose art made hyperbolic phenomena familiar 
to the general public. On the other hand there are many theoretical results on Kleinian groups, the symmetries hyperbolic 3-space, but only very few illustrations. 

The boundary of hyperbolic $3$-space $\partial \mathbb{H}^3$ may be identified with the Riemann Sphere 
$\widehat{\mathbb{C}}=\mathbb{C}\cup\{\infty\}$, or alternatively the complex projective plane $\textbf{P}({\mathbb{C}^1})$. The conformal 
automorphisms of $\widehat{\mathbb{C}}$ onto itself are given by M\"obius transformations $z \mapsto \frac{az+b}{cz+d}$ where 
$a,b,c,d \in \mathbb{C}$ with $ad-bc\neq 0$. 
If we consider the conformal ball model (Poincar\'e model) of hyperbolic space, the M\"obius maps naturally extend from the boundary of the space to 
isometries of hyperbolic 3-space. Kleinian groups are discrete subgroups of M\"obius transformations, a special class of which 
are the Coxeter groups of hyperbolic 3-space.

Kleinian groups have fascinated the imagination of many mathematicians, who have created many beautiful images on the complex plane and 
Reimann sphere. 
Perhaps the first such fractal image is due to Felix Klein himself, who together with Fricke in the late 19th century hand illustrated the 
limit set of a Schottky group featuring the first few iterations of an infinite chain of tangent circles \cite{KF1897}. 

Later during the advent of computer graphics such fractal images were popularized by Benoit Mandelbrot resulting in great renewed interest in the field \cite{M}. David Mumford's book {\it Indra's Pearls} gives an outstanding introduction to the subject for a boarder mathematically inclined audience with ample figures, complete with instructions on how to reproduce them \cite{Mum}.

The reader will find that our images of optimal horoball packings resemble such fractal images. 
This should not be surprising as their group of symmetries are isomorphic to Kleinian groups. We shall use the projective Beltrami-Klein model in 
order to preserve the convexity of the polyhedral cells of our tilings. Such tilings have as their group of isometries discrete 
subgroups of $\mathbf{PGL}(4,\mathbb{R})$. By a suitable isometry the Beltrami-Klein model and the conformal model are equivalent.

In this paper we illustrate our results found in \cite{KSz}, in which we give the necessary mathematical background by proving the existence of multiple optimal packing configurations in extended hyperbolic space $\overline{\mathbb{H}}^3$. We describe a procedure to study and visualize the densest horoball packing arrangements using their group of symmetries. 
These horoball packings are related to the regular Coxeter honeycombs with Schl\"afli symbols $\{3,3,6 \}$ and $\{4,3,6\}$  
where we allow horoballs of different types, and the resulting images resemble Apollonian gaskets.
Our technique is suitable to describe other hyperbolic ball or horoball packings generated by arbitrary hyperbolic Coxeter groups described in \cite{C54}.  
We implement our method in Mathematica, the Wolfram Language. The code used to generate the images is freely accessible from the personal webpage of the authors \cite{K}. 

\section{Background and Motivation}

In hyperbolic $3$-space the densest possible ball packings are realized as horoball packings \cite{BF64, KSz}. In $\mathbb{H}^3$ a horosphere is the higher dimensional analog of a horocycle in $\mathbb{H}^2$, in a sense that will be made precise in Section 2. A horoball is a horosphere together with its interior, and a horoball packing $\cB$ of $\overline{\mathbb{H}}^3$ is a countable collection of non-overlapping horoballs in $\overline{\mathbb{H}}^3$, i.e. for any $B_1,B_2 \in \cB$ it holds that $\text{Int}(B_1 \cap B_2)=\emptyset$.

The definition of packing density is critical in hyperbolic space as shown by the famous example of B\"or\"oczky \cite{B78, GK}. 
The most widely accepted notion of packing density, considers local densities of balls with respect to their Dirichlet--Voronoi cells 
(cf. \cite{B78} and \cite{K98}). 
Let $B \in \cB$ be a horoball of the packing, and $p \in \overline{\mathbb{H}}^3$ be an arbitrary point.
Define $d(p,B)$ to be the perpendicular distance from point $p$ to the horosphere $S = \partial B$, where $d(p,B)$
is negative for $p \in B$. The Dirichlet--Voronoi cell $\cD(B,\cB)$ of horoball $B$ with respect to packing $\cB$ is defined as the convex body
\begin{equation}
\cD(B,\cB) = \{ p \in \mathbb{H}^3 | d(p,B) \le d(p,B'), ~ \forall B' \in \cB \}.
\notag
\end{equation}
Both $B$ and $\cD$ may have unbounded volume, so the usual notion of local density is
modified as follows. Let $q \in \partial{\mathbb{H}}^3$ denote the ideal center of $B$, and take 
its boundary $S$ to be the one-point compactification of Euclidean plane.
Let $D_c(r) \subset S$ be a disk with center $c \in S \setminus \{q\}$.
Then $q \in \partial {\mathbb{H}^3}$ and $D_c(r)$ determine a convex cone
$\cC(r) = cone_q(D_c(r)) \in \overline{\mathbb{H}}^3$ with
apex $q$ consisting of all hyperbolic geodesics passing through $D_c(r)$ with limit point $q$. 
The local density $\delta(B, \cB)$ of $B$ to $\cD$ is defined as
\begin{equation}
\delta(\cB, B) =\varlimsup\limits_{r \rightarrow \infty} \frac{vol(B \cap \cC(r))} {vol(\cD \cap \cC(r))}.
\notag
\end{equation}
This limit is independent of the choice of center for $D_c(r)$.

In the case of periodic ball or horoball packings, the local density defined above can be extended to the entire hyperbolic space. 
This local density
is related to the simplicial density function \cite{K98}, generalized in \cite{Sz13} and \cite{Sz12-2}.
In this paper we will use such definition of packing density (cf. Section 3).  

Horoballs are always congruent in $\overline{\mathbb{H}}^3$ in the classical sense of isometries. 
In \cite{Sz12-2} we refined the notion of ``congruent" horoballs. 
Two horoballs of a packing are of the ``same type", or ``equipacked", if and only if the
local densities of each horoball with respect to the ambient cell (e.g. D-V cell, or fundamental domain) are equal,
otherwise they are said to be of ``different type".

If we assume horoballs belong to the ``same type" then by analytical continuation
the simplicial density function on $\overline{\mathbb{H}}^n$ $(n \geq 2)$ can be extended from balls of radius $r < \infty$ to 
the infinite case.  In particular consider $n + 1$ horoballs which are mutually tangent, then the convex hull of their
base points at infinity will be a totally asymptotic (or ideal) regular simplex $T_{reg}^{\infty} \in \overline{\mathbb{H}}^n$ of finite volume.
Let $B$ be one such horoball, then
\begin{equation}
d_n(\infty) = (n + 1)\frac{vol(B \cap T_{reg}^\infty)}{vol(T_{reg}^\infty)}.
\notag
\end{equation}
For a horoball packing $\cB$, there is an analogue of ball packing, namely (cf. \cite{B78}, Theorem 4)
\begin{equation}
\delta_n(\cB, B) \le d_n(\infty),~ \forall B \in \cB.
\notag
\end{equation}
The upper bound $d_n(\infty)$ $(n=2,3)$ is attained for a regular horoball packing, that is, a
packing by horoballs which are inscribed in the cells of a regular honeycomb of $\overline{\mathbb{H}}^n$. For
dimensions $n = 2$, there is only one such packing. It belongs to the regular tessellation $\{\infty, 3 \}$ . Its dual
$\{3,\infty\}$ is the regular tessellation by ideal triangles all of whose vertices are surrounded
by infinitely many triangles. This packing has in-circle density $d_2(\infty)=\frac{3}{\pi} \approx 0.95493 $.

In $\overline{\mathbb{H}}^3$ there is exactly one horoball packing with horoballs of the same type whose Dirichlet--Voronoi cells give rise to a
regular honeycomb described by Schl\"afli symbol $\{6,3,3\}$ . Its
dual $\{3,3,6\}$ consists of ideal regular simplices $T_{reg}^\infty$  with dihedral angles $\frac{\pi}{3}$ that form a 6-cycle around each edge
of the tessellation. The density of this packing is $\delta_3 (\infty) \approx 0.85328$.

If horoballs of different types are allowed at the ideal vertices, then we generalize
the notion of the simplicial density function \cite{Sz12-2}.
In \cite{KSz}  we gave several new examples of horoball packing arrangements based on totally asymptotic Coxeter tilings that yield the
B\"or\"oczky--Florian upper bound \cite{BF64}, showing that the optimal ball packing arrangement in $\mathbb{H}^3$ described above is not unique.

In \cite{KSz14} we investigated ball packings in hyperbolic $4$-space. Using the techniques described
we found several counterexamples to a conjecture of L. Fejes-T\'oth regarding the upper bound of packing density in $\mathbb{H}^4$ \cite{FTL} . The highest known packing density now is $\approx 0.71645$. The hyperbolic regular $24$-cell and its
regular $4$-dimensional honeycomb with Schl\"afli symbol $\{3,4,3,4\}$ also yields this new optimal packing density.

In addition, in \cite{Sz12-2, Sz13} we found that
by admitting horoballs of different types at each vertex of a totally asymptotic simplex, we locally exceed the B\"or\"oczky-type density upper bound for $n > 3$.
For example, in $\overline{\mathbb{H}}^4$ the locally optimal packing density was found to be
$\approx 0.77038$, higher than the B\"or\"oczky-type density upper bound $\approx 0.73046$.
However such packings are only locally optimal and cannot be extended to pack the entire $\mathbb{H}^n$.
%
\subsection{Projective Geometry of $\mathbb{H}^3$}
%
In what follows we use the Beltrami-Klein model, and a projective interpretation of hyperbolic geometry. As we are primarily interested 
in packings of tilings with convex polyhedral cells, this model has the advantage of greatly 
simplifying our density calculations compared to conformal models such as the Poincar\'e model where convexity is 
severely distorted \cite{KSz}. 
In this section we give a brief review of the concepts used in this paper. For a general discussion and 
background in hyperbolic geometry, as well as the projective models of the eight Thurston geometries see \cite{M97}. For a general higher 
dimensional discussion and examples in hyperbolic 4-space see \cite{KSz14}.

\subsection{The Projective Model}
We use the projective model in Lorentzian $4$-space
$\mathbb{E}^{1,3}$ of signature $(1,3)$, i.e.~$\mathbb{E}^{1,3}$ is
the real vector space $\mathbf{V}^{4}$ equipped with the bilinear
form of signature $(1,3)$
\begin{equation}
\langle  \mathbf{x},\mathbf{y} \rangle = -x^0y^0+x^1y^1 +x^2y^2 + x^3 y^3 \label{bilinear_form} \tag{2.1}
\end{equation}
where the non-zero real vectors
$\mathbf{x}=(x^0,x^1,x^2,x^3)\in\mathbf{V}^{4}$
and $ \mathbf{y}=(y^0,y^1,y^2, y^3)\in\mathbf{V}^{4}$ represent points in projective space
$\mathbf{P}^3(\mathbb{R})$. $\mathbb{H}^3$ is represented as the
interior of the absolute quadratic form
\begin{equation}
Q=\{[\mathbf{x}]\in\mathbf{P}^3 | \langle \mathbf{x},\mathbf{x} \rangle =0 \}=\partial \mathbb{H}^3  \tag{2.2}
\end{equation}
in real projective space $\mathbf{P}^3(\mathbf{V}^{4},\mbox{\boldmath$V$}\!_{4})$. 
All proper interior points $\mathbf{x} \in \mathbb{H}^3$ are characterized by
$\langle \mathbf{x},\mathbf{x} \rangle < 0$.

The boundary points $\partial \mathbb{H}^3 $ in
$\mathbf{P}^3$ represent the absolute points at infinity of $\mathbb{H}^3$.
Points $\mathbf{y}$ satisfying $\langle \mathbf{y}, \mathbf{y} \rangle >
0$ lie outside $\partial \mathbb{H}^3$ and are called the outer points
of $\mathbb{H}^3$. Take $P([\mathbf{x}]) \in \mathbf{P}^3$, point
$[\mathbf{y}] \in \mathbf{P}^3$ is said to be conjugate to
$[\mathbf{x}]$ relative to $Q$ when $\langle
\mathbf{x}, \mathbf{y} \rangle =0$. The set of all points conjugate
to $P([\mathbf{x}])$ form a projective (polar) hyperplane
\begin{equation}
pol(P):=\{[\mathbf{y}]\in\mathbf{P}^3 | \langle \mathbf{x},\mathbf{y} \rangle =0 \}. \tag{2.3}
\end{equation}
Hence the bilinear form $Q$ in (\ref{bilinear_form}) induces a bijection
or linear polarity $\mathbf{V}^{4} \rightarrow
\mbox{\boldmath$V$}\!_{4}$
between the points of $\mathbf{P}^3$
and its hyperplane.
Point $X [\bold{x}]$ and hyperplane $\alpha
[\mbox{\boldmath$a$}]$ are incident if the value of
linear form $\mbox{\boldmath$a$}$ evaluated on vector $\bold{x}$ is
 zero, i.e. $\bold{x}\mbox{\boldmath$a$}=0$ where $\mathbf{x} \in \
\mathbf{V}^{4} \setminus \{\mathbf{0}\}$, and $\ \mbox{\boldmath$a$} \in
\mbox{\boldmath$V$}_{4} \setminus \{\mbox{\boldmath$0$}\}$.
Similarly, lines in $\mathbf{P}^3$ are characterized by
2-planes of $\mathbf{V}^{4}$ or $2$-planes of $\
\mbox{\boldmath$V$}\!_{4}$ \cite{M97}.

Let $P \subset \mathbb{H}^3$ denote a polyhedron bounded by
a finite set of hyperplanes $H^i$ with unit normal vectors
$\mbox{\boldmath$b$}^i \in \mbox{\boldmath$V$}\!_{n+1}$ directed
 towards the interior of $P$:
\begin{equation}
H^i:=\{\mathbf{x} \in \mathbb{H}^3 | \langle  \mathbf{x},\mbox{\boldmath$b$}^i \rangle =0 \} \ \ \text{with} \ \
\langle \mbox{\boldmath$b$}^i,\mbox{\boldmath$b$}^i \rangle = 1. \tag{2.4}
\end{equation}
In this paper $P$ is assumed to be an acute-angled polyhedron
with proper or ideal vertices.
The Grammian matrix $G(P):=( \langle \mbox{\boldmath$b$}^i,
\mbox{\boldmath$b$}^j \rangle )_{i,j} ~ {i,j \in \{ 0,1,2,3 \} }$  is an
indecomposable symmetric matrix of signature $(1,3)$ with entries
$\langle \mbox{\boldmath$b$}^i,\mbox{\boldmath$b$}^i \rangle = 1$
and $\langle \mbox{\boldmath$b$}^i,\mbox{\boldmath$b$}^j \rangle
\leq 0$ for $i \ne j$ where
$$
\langle \mbox{\boldmath$b$}^i,\mbox{\boldmath$b$}^j \rangle =
\left\{
\begin{aligned}
&0 & &\text{if}~H^i \perp H^j,\\
&-\cos{\alpha^{ij}} & &\text{if}~H^i,H^j ~ \text{intersect \ along an edge of $P$ \ at \ angle} \ \alpha^{ij}, \\
&-1 & &\text{if}~\ H^i,H^j ~ \text{are parallel in the hyperbolic sense}, \\
&-\cosh{l^{ij}} & &\text{if}~H^i,H^j ~ \text{admit a common perpendicular of length} \ l^{ij}.
\end{aligned}
\right.
$$
This information is encoded in the weighted graph or scheme of the polytope $\sum(P)$. The graph nodes correspond
to the hyperplanes $H^i$ and are connected if $H^i$ and $H^j$ not perpendicular ($i \neq j$).
If they are connected we write the positive weight $k$ where  $\alpha_{ij} = \pi / k$ on the edge, and
unlabeled edges denote an angle of $\pi/3$. This graph is also known as the Coxeter--Dynkin diagram

In this paper we set the sectional curvature of $\mathbb{H}^3$,
$K=-k^2$, to be $k=1$. The distance $d$ of two proper points
$[\mathbf{x}]$ and $[\mathbf{y}]$ is given by
\begin{equation}
\cosh \left(d ([\mathbf{x}],[\mathbf{y}]) \right)=\frac{-\langle ~ \mathbf{x},~\mathbf{y} \rangle }{\sqrt{\langle ~ \mathbf{x},~\mathbf{x} \rangle
\langle ~ \mathbf{y},~\mathbf{y} \rangle }}. \tag{2.5}
\end{equation}
The perpendicular foot $Y[\mathbf{y}]$ of point $X[\mathbf{x}]$ dropped onto plane $[\mbox{\boldmath$u$}]$ is given by
\begin{equation}
\mathbf{y} = \mathbf{Pr_{u}(x)} = \mathbf{x} - \frac{\langle \mathbf{x}, \mathbf{u} \rangle}{\langle \mathbf{u}, \mathbf{u} \rangle} \mathbf{u}, \tag{2.6}
\end{equation}
where $[\mathbf{u}]$ is the pole of the plane $[\mbox{\boldmath$u$}]$.

\subsection{Characterization of horoballs in $\mathbb{H}^3$}

A horosphere in $\mathbb{H}^3$ is a
hyperbolic $2$-sphere with infinite radius that is centered
at an ideal point, on $\partial \mathbb{H}^3$. Equivalently, a horosphere is a $2$-surface orthogonal to
the set of parallel geodesics passing through a point of the absolute quadratic surface.
A horoball is a horosphere together with its interior.

We consider the usual Beltrami-Klein ball model of $\mathbb{H}^3$
centered at $O(1,0,0,0)$ with a given vector basis
$\bold{a}_i$ where $i=0,1,2, 3$ and set an
arbitrary point at infinity to lie at $T_0=(1,0,0,1)$.
The equation of a horosphere with center
$T_0=(1,0,0,1)$ passing through point $S=(1,0,0,s)$ is derived from the equation of the
the absolute sphere $||\mathbf{x}||^2 = -x^0 x^0 +x^1 x^1+x^2 x^2+ x^3 x^3 = 0$, and the plane $x^0-x^3=0$ tangent to the absolute sphere at $T_0$.
The general equation of the horosphere in projective coordinates is
\begin{align}
(s-1)\left( -x^0 x^0 +x^1 x^1+x^2 x^2+ x^3 x^3\right)-(1+s){(x^0-x^3)}^2 & =0, \tag{2.7}
\end{align}

where $s \neq \pm1$. The equation for the horophere in Cartesian coordinates is obtained by setting 
$x=\frac{x^1}{x^0}$, $y=\frac{x^2}{x^0}$, and $z=\frac{x^3}{x^0}$,
\begin{equation}
\frac{2(x^2+y^2)}{1-s}+\frac{4(z-\frac{s+1}{2})^2}{{(1-s)}^2}=1. \tag{2.8}
\end{equation}

For polar plots it is useful to have the polar form of the horosphere equation with parameters $s \in (-1,1)$,
$\phi\in[0, 2\pi)$, and $\theta\in[0,\pi]$,
\begin{equation}
\begin{gathered}
x=\sqrt{\frac{1-s}{2}}\sin{\theta} \cos{\phi}, ~~~~
y=\sqrt{\frac{1-s}{2}}\sin{\theta} \sin{\phi}, \\
z=\frac{1+s}{2}+\frac{1-s}{2}\cos\theta .
\end{gathered} \tag{2.9}
\end{equation}

Applying rotations to these equations one can obtain the equations of horospheres centered at
an arbitrary point on the boundary of the model.

In $\mathbb{H}^3$ any two horoballs are congruent in the classical sense, there exists a hyperbolic isometry mapping one to another.
However, in our approach we find it rewarding to distinguish between certain horoballs of a packing.
We shall use the notion of horoball type with respect to a packing as introduced in \cite{Sz12-2}. The motivation is that one has a one-parameter family of concentric horoballs centered at each ideal point of the boundary of the model sphere. Indeed, each horoball in such family corresponds to a different value of parameter $s\in(-1,1)$ in the above equations. Concentric horoballs with different $s$-parameters may have different relative densities with respect to the fundamental domain of the packing.

\begin{definition}
Two horoballs of a regular horoball packing are of the {\it same type} or {\it equipacked} if
and only if their local packing densities with respect the fundamental domain are equal.
Otherwise the horoballs are of {\it different type}.
\end{definition}

The hyperbolic length $l(x)$ of a horospheric arc belonging to a chord segment of length $x$ is given by
$l(x) = 2 \sinh{\left(\frac{x}{2}\right)}$.
The  intrinsic geometry of a horosphere is Euclidean,
so the $2$-dimensional volume $\mathcal{A}$ of a region $A$ on the
surface of a horosphere is calculated as in $\mathbb{E}^{2}$.
The volume of the horoball piece $\mathcal{H}(A)$ determined by $A$ and
the aggregate of axes
drawn from $A$ to the center of the horoball is
\begin{equation}
\label{eq:bolyai}
Vol(\mathcal{H}(A)) = \frac{1}{2}\mathcal{A}. \tag{2.10}
\end{equation}

\section{Visualization of the Optimal Packings}

A regular packing is fully determined by the ball arrangement in its fundamental 
domain. 
Our method for visualization of horoball arrangements is based on the use of Coxeter groups which are the symmetries of our packings.
A fundamental domain of the Coxeter group is an orthoscheme of degree $0$ with given dihedral angles.

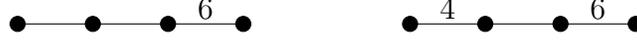
\begin{figure}[p]
\begin{center}
   \begin{tikzpicture}

	\draw (0,0) -- (3,0);
	
	\draw[fill=black] (0,0) circle(.1);
	\draw[fill=black] (1,0) circle(.1);
	\draw[fill=black] (2,0) circle(.1);
	\draw[fill=black] (3,0) circle(.1);
	
	\node at (2.5,0.2) {$6$};
	 
\end{tikzpicture}~~~~~~~~~~~~~~~~~~
\begin{tikzpicture}

	\draw (0,0) -- (3,0);
	
	\draw[fill=black] (0,0) circle(.1);
	\draw[fill=black] (1,0) circle(.1);
	\draw[fill=black] (2,0) circle(.1);
	\draw[fill=black] (3,0) circle(.1);

	\node at (.5,0.2) {$4$};	
	\node at (2.5,0.2) {$6$};
	 
\end{tikzpicture}
\caption{Coxeter--Dynkin diagrams of $\{3,3,6\}$ and $\{4,3,6\}$.}
\label{Cox-Dyn}
\end{center}
\end{figure}

In the case of the four optimal horoball packings in $\mathbb{H}^3$, the centers of the horoballs are arranged at the lattice points of $\mathbb{H}^3$ tiled by either totally asymptotic regular hyperbolic 
tetrahedra with dihedral angles $\frac{\pi}{3}$, or totally asymptotic regular hyperbolic cubes with dihedral angles $\frac{\pi}{5}$. 
Therefore, in order to find the centers of the horoballs of the packings 
we use the elements $\bg_i$ of the generator set of the Coxeter group of the tilings. The ideal regular tetrahedron or regular cube are the fundamental domains of the Coxeter groups of the tilings that preserve our packings. The subgroup corresponding to the tetrahedron is called tetrahedral group of isometries and is denoted by $\bT$.
The subgroup that belongs to the cube is called cubic group of isometries and is denoted by $\bC$. For the schemes of the tilings see Fig. \ref{Cox-Dyn}.

The Coxeter group is then used to iteratively generate the packing by successively applying generators to map the packing configuration within the fundamental domain to all of $\mathbb{H}^3$.
The metric data we use to describe the four packings in this discussion is consistent with that used in \cite{KSz}, where we proved the optimality of these packings. We shall assume all statements of optimality given in \cite{KSz}, and omit any proofs. 

\begin{figure}[p]
\begin{center}
\caption{The two optimal horoball packings of the tiling $\{4,3,6\}$, with packings generated by 3, 4, and 5 reflections (crowns, or layers about the base cube) of horoballs about the fundamental domain.}
\begin{tabular}{|c|c|c|}
\hline
Crowns & Maximal Case & Balanced Case \\
\hline
3&
\includegraphics[height=60mm]{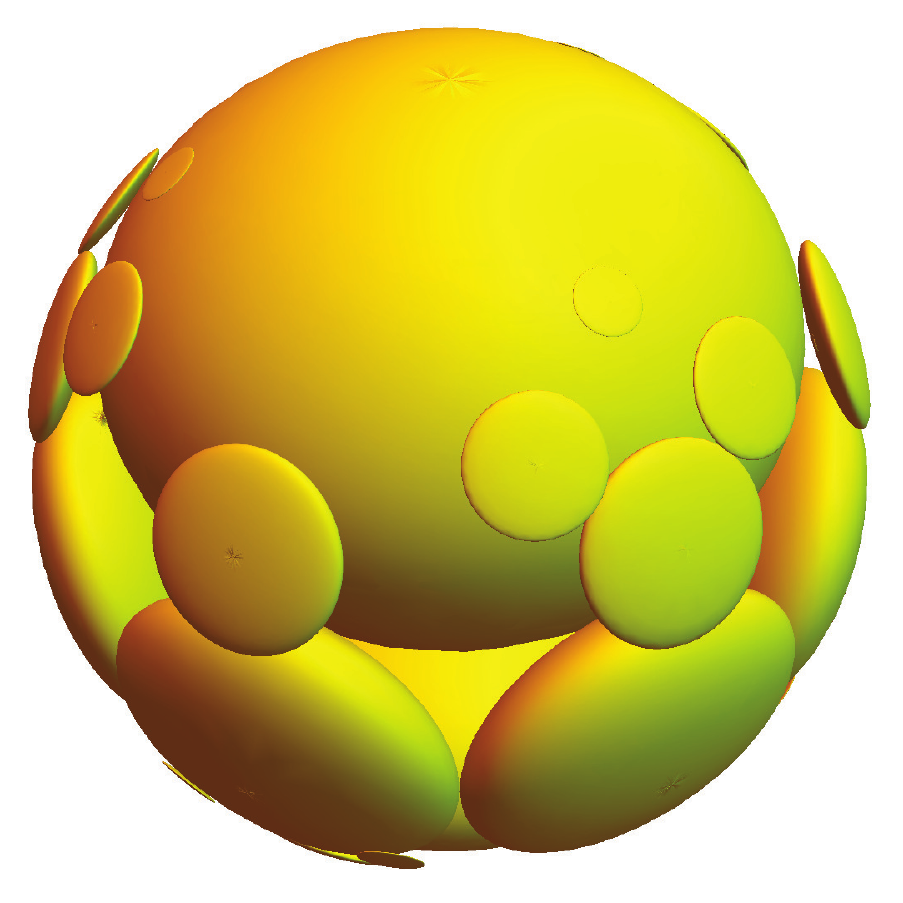} &
\includegraphics[height=60mm]{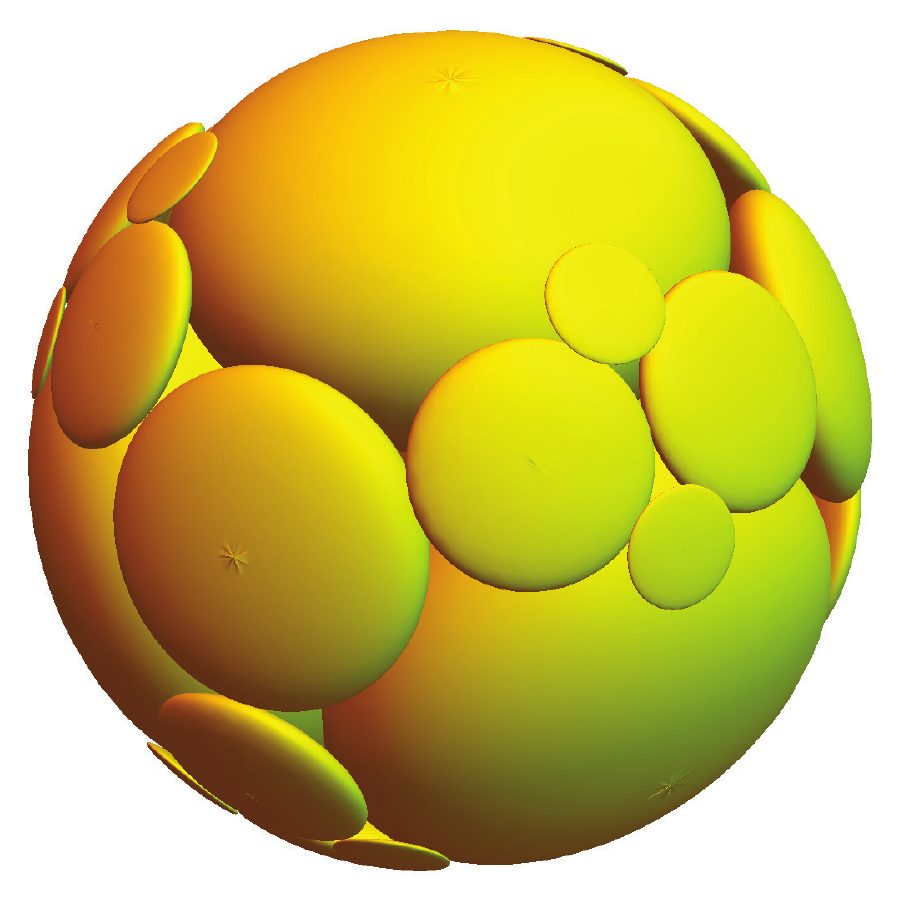}\\
\hline
4&
\includegraphics[height=60mm]{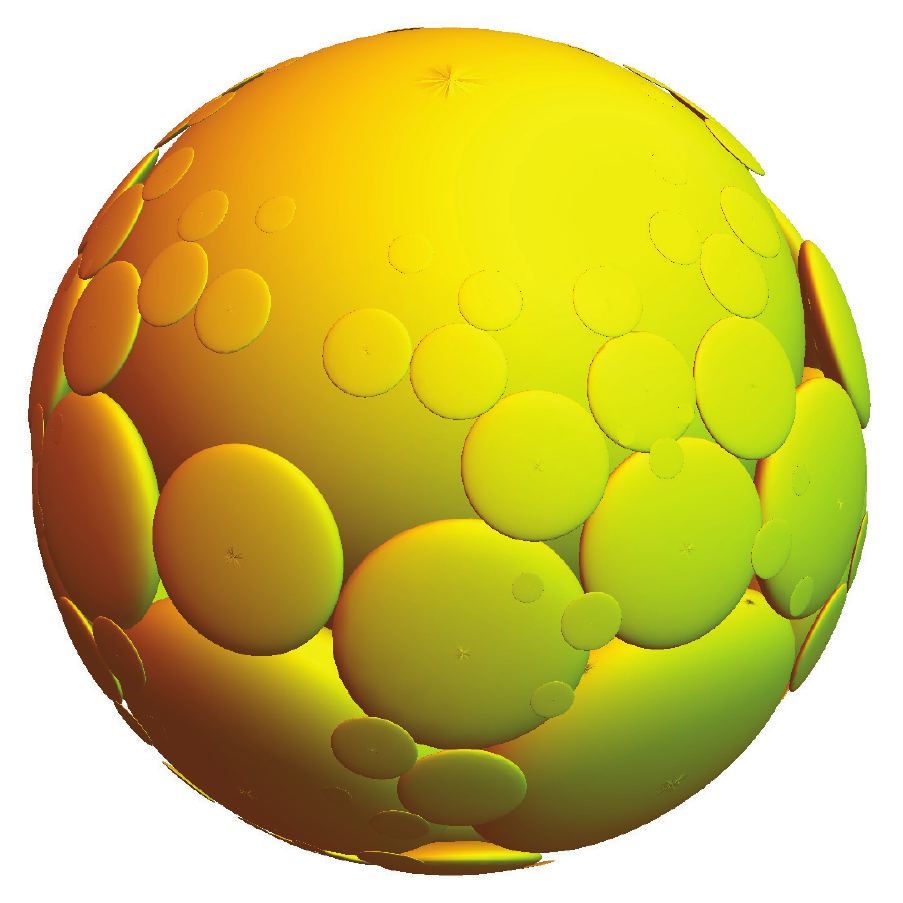} &
\includegraphics[height=60mm]{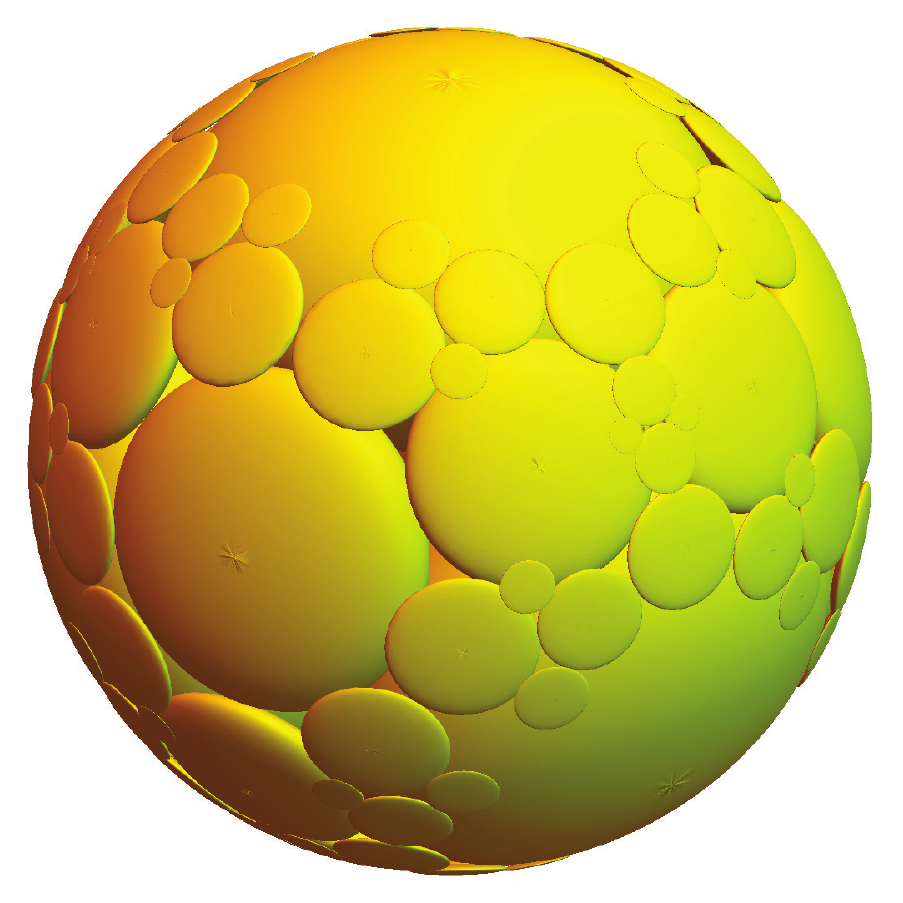}\\
\hline
5&
\includegraphics[height=60mm]{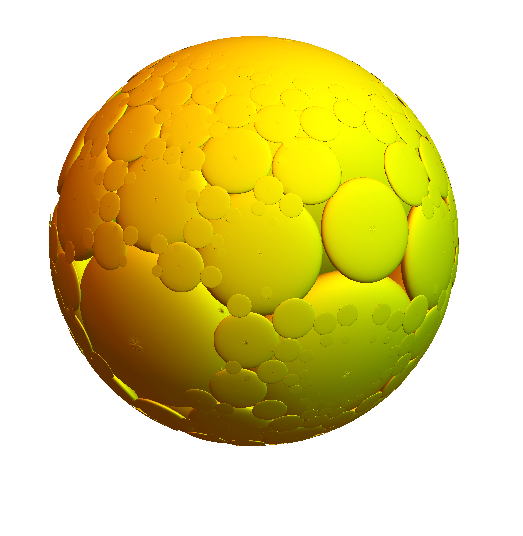}&
\includegraphics[height=60mm]{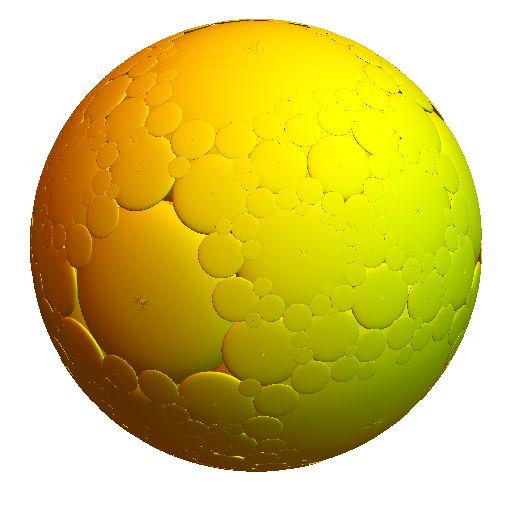}\\
\hline
\end{tabular}

\label{fig:436}
\end{center}
\end{figure}

For the Mathematica code used to generate Figures \ref{fig:336} and \ref{fig:436} see \cite{K}.

\subsection{$\{3, 3, 6\}$ Tetrahedral Tiling}

The $\{3, 3, 6\}$ Coxeter tiling is a three dimensional honeycomb with cells consisting of fully asymptotic
regular tetrahedra. The two extremal cases of horoball configurations yield the optimally dense packings (see \cite{KSz}).  
In this section we will restrict our attention to these two cases.

\subsubsection{Fundamental Domain}

To parameterize the fundamental domain of the subtiling of the Coxeter honeycomb, fix a regular totally asymptotic tetrahedron 
$E_0E_1E_2E_3$ as the fundamental domain. Horoballs are centered at vertices $E_0, \dots,
E_3$ so that they preserve symmetries of the packing preserve the fundamental domain. 
The two optimal packing configurations for this case were found in \cite{KSz}.

The barycentric subdivision of one tetrahedral cell gives six congruent orthoschemes.
Define orthoscheme $A_1A_2A_3A_4$ by setting $A_1=E_0$,
$A_4=E_1$, $A_3$ the center of the triangular facet $E_1E_2E_3$ opposite vertex $E_0$, 
and $A_2$ to be the perpendicular foot of $E_0$ projected onto edge $E_1E_2$.  The
Schl\"afli symbol of orthoscheme $A_1A_2A_3A_4$ is $\{3,6,3\}$.
A metric description of the fundamental domain is given by assigning coordinates 
\begin{equation*}
\begin{gathered}
A_1=(1,0,0,1), ~ A_2 =\left(1,\frac{\sqrt{3}}{4},\frac{1}{4},0\right),~ A_3=(1,0,0,0),~A_4=(1,0,1,0),
\end{gathered}
\end{equation*}
to orthoscheme $A_0A_1A_2A_3$. The associated tetrahedral cell then has coordinates
{\small
\begin{equation*}
\begin{gathered}
E_0 =(1, 0, 0, 1),~E_1=(1, 0, 1, 0),~E_2 =\left( 1,\frac{\sqrt{3}}{2},-\frac{1}{2},0\right),
~E_3 =\left( 1,-\frac{\sqrt{3}}{2},-\frac{1}{2},0\right),
\end{gathered}
\end{equation*}}
which give the fundamental group of the packing.
One may check the angle requirements of the tiling are satisfied by computing the inner products of the normals $[\mathbf{u}_i]$ in Table \ref{table:data_336}.

A Coxeter group $\bG$ is a finitely generated group defined by a presentation of the form $\langle~\bg_i~ | ~(\bg_i~\bg_j)^{k_{ij}} \rangle$ where $k_{ij}$ is a positive integer or $\infty$ satisfying well known symmetry assumptions.
The Coxeter group $\bT$ acts by isometries (or congruence transformations) of $\mathbb{H}^3$.
The four generators $\{\bg_i\}_{i=1}^4$ of $\bT$ are determined by reflections on the four sides of the fundamental domain of the packing. The 
reflecting planes themselves are uniquely determined by the choice of coordinates for the regular tetrahedron. 

The symmetries of the 
Beltrami-Klein model are given by $\mathbf{PSL}(4,\mathbb{R})$ so it remains to find matrix representations for the $\bg_i$. All vertices of the fundamental domain are ideal, hence all lattice points generated by the group are also ideal. 
In particular, let $\bg_i$ be the reflection across the plane of the facet opposite vertex $E_i$.
Then $\bg_i(E_i)$ is the intersection with the model sphere of the geodesic passing through the two points $E_i$ and  $\mathbf{Pr}_{\bu_i}(E_i)$, 
the perpendicular foot of the vertex projected onto the facet with normal $[\bu_i]$. We compute the results of the actions to find $\bg_i(E_i)$ as in Table \ref{table:data_336}. 
Here $\bg_i \in \mathbf{PSL}(4,\mathbb{R})$ is the group generator corresponding to the $i$-th basis element. We next find the matrix form of the reflections $\bg_i$ by using the set 
of eigenvectors consistent with reflection onto facet opposite $E_i$. Such a matrix leaves the plane of the facet invariant, so to find $\bg_i$ we set $\bv_i = [E_i]$, and compute the solutions to the linear system

\begin{equation}
\bg_i.\bv_j=\lambda_j \bw_j 
\end{equation} 
where $\bw_j=\bv_j$ if $i \neq j$ and $\bg_i(E_i)$ otherwise. 
The data used for these computations is summarized in Table \ref{table:data_336}.

The set of generators of the subgroup $\bT$ as 
elements of $\mathbf{PSL}(4,\mathbb{R})$ consistent with choice of coordinates $E_i$ is 
\begin{equation*}
\begin{gathered}
\bg_1
=\left(
\begin{array}{cccc}
 1 & 0 & 0 & 0 \\
 0 & 1 & 0 & 0 \\
 0 & 0 & 1 & 0 \\
 0 & 0 & 0 & -1 \\
\end{array}
\right),
~
\bg_2=\left(
\begin{array}{cccc}
 -\frac{3}{2} & 0 & -1 & \frac{1}{2} \\
 0 & -1 & 0 & 0 \\
 1 & 0 & 1 & -1 \\
 -\frac{1}{2} & 0 & -1 & -\frac{1}{2} \\
\end{array}
\right),
\\
\bg_3=\left(
\begin{array}{cccc}
 -\frac{3}{2} & -\frac{\sqrt{3}}{2} & \frac{1}{2} & \frac{1}{2} \\
 \frac{\sqrt{3}}{2} & \frac{1}{2} & -\frac{\sqrt{3}}{2} & -\frac{\sqrt{3}}{2} \\
 -\frac{1}{2} & -\frac{\sqrt{3}}{2} & -\frac{1}{2} & \frac{1}{2} \\
 -\frac{1}{2} & -\frac{\sqrt{3}}{2} & \frac{1}{2} & -\frac{1}{2} \\
\end{array}
\right),
~
\bg_4=\left(
\begin{array}{cccc}
 -\frac{3}{2} & \frac{\sqrt{3}}{2} & \frac{1}{2} & \frac{1}{2} \\
 -\frac{\sqrt{3}}{2} & \frac{1}{2} & \frac{\sqrt{3}}{2} & \frac{\sqrt{3}}{2} \\
 -\frac{1}{2} & \frac{\sqrt{3}}{2} & -\frac{1}{2} & \frac{1}{2} \\
 -\frac{1}{2} & \frac{\sqrt{3}}{2} & \frac{1}{2} & -\frac{1}{2} \\
\end{array}
\right).
\end{gathered}
\end{equation*}
\subsubsection{Horoball Packings of the Fundamental Domain}

\begin{figure}[p]
\begin{center}
\caption{Optimal horoball configurations in the fundamental domain of the $\{3,3,6\}$ tiling. The outer regular asymptotic tetrahedron represents the fundamental domain. The barycentric simplex is the orthoscheme $\{6,3,6\}$ of the tiling. Horoballs are centered at each vertex. The outer sphere represents the boundary of $\partial \mathbb{H}^3$ in the Beltrami--Klein model.}
\begin{tabular}{cc}
\includegraphics[height=60mm]{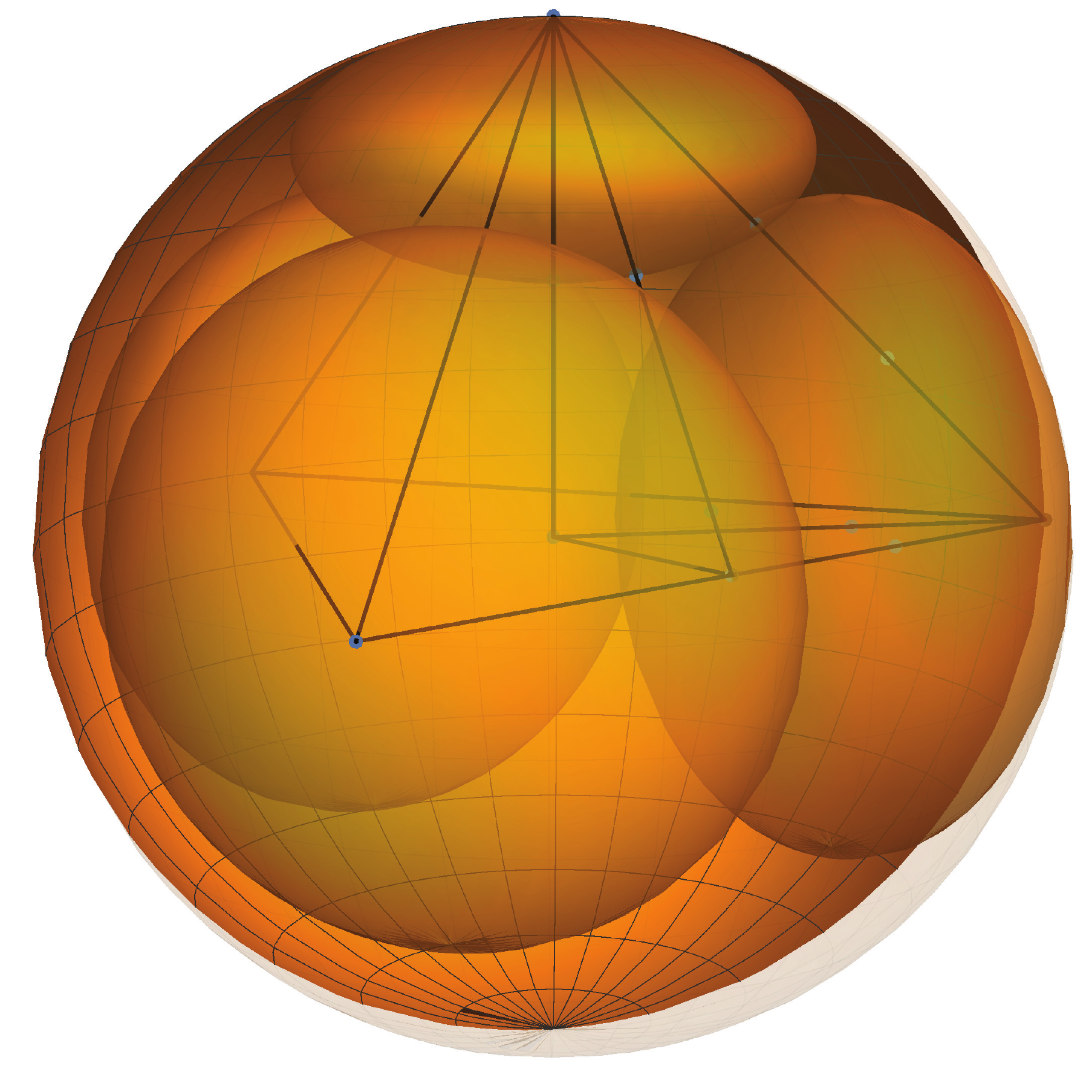} &
\includegraphics[height=60mm]{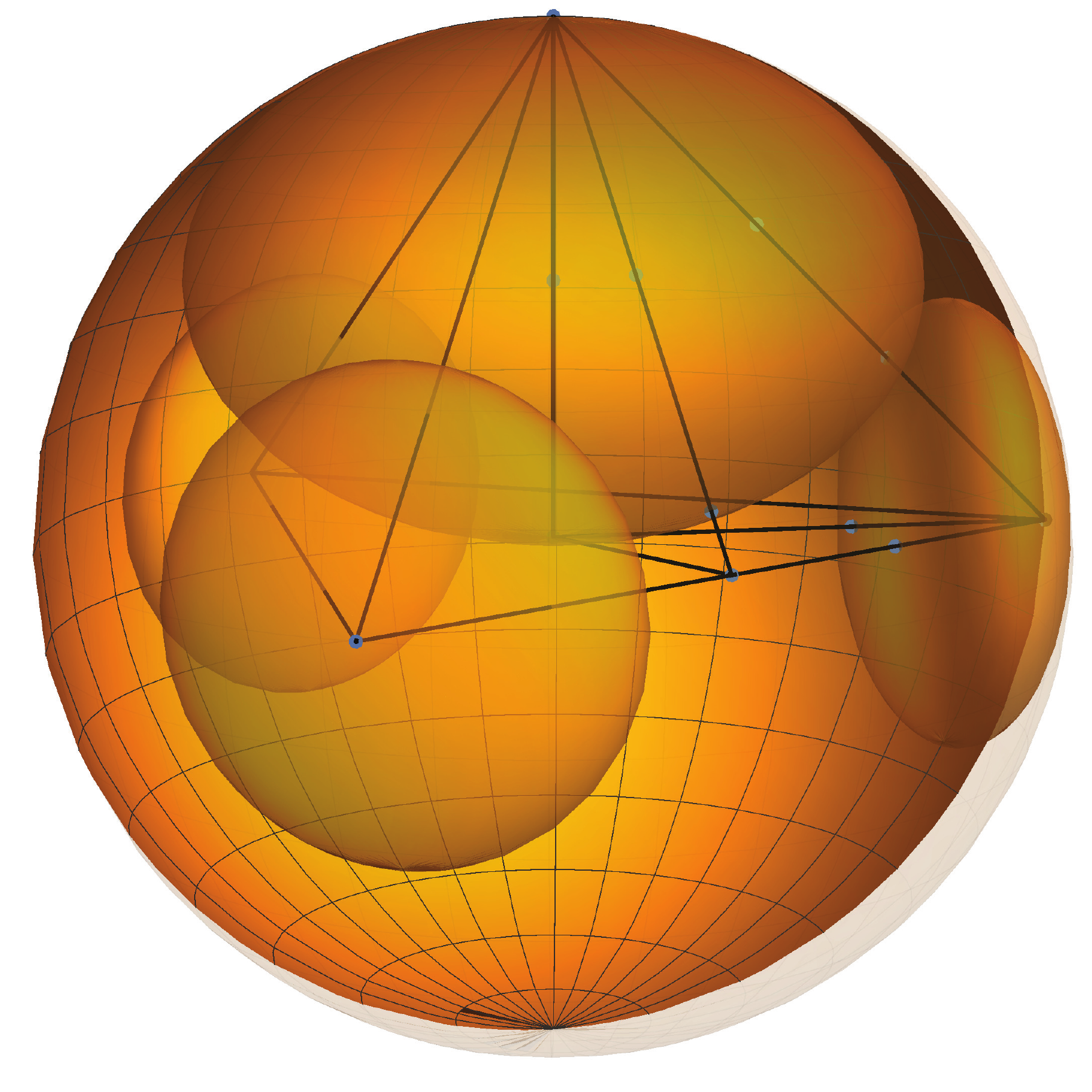} \\
(a) B\"or\"oczky--Florian case. & (b)  Kozma--Szirmai case.\\
\end{tabular}
\label{fig:336model}
\end{center}
\end{figure}

Let $B_0$ and $B_3$ be two horoballs centered at $E_0$ and $E_3$, i.e.,
the two vertices of the tetrahedra common with the orthoscheme. The
density of the $\{3,3,6\}$ Coxeter tiling is defined by 
\begin{equation}
\delta(\mathcal{B}_{336})=\frac{Vol(B_0\cap\mathcal{O}_{(3,6,3)})+Vol(B_3\cap\mathcal{O}_{(3,6,3)})}{Vol({\mathcal{O}_{(3,6,3)}})}.
\end{equation}
\begin{proposition}
The packing density obtained in orthoscheme $\mathcal{O}_{(3,6,3)}$ can be
extended to tetrahedron $\bT$ and therefore to the entire
$\mathbb{H}^3$.
\end{proposition}
There are two cases yielding the optimal packing density of $0.85328\dots$:
\begin{enumerate}
\item B\"or\"oczky--Florian Case: This represents the equilibrium case where all horoballs are equipacked with respect to the fundamental 
domain. Horoballs meet along the midpoint of each edge in the model sphere. The data for the horoballs is using the s parameter 
$s_0=0$ and $s_1=\frac{3}{5}$. See Fig. \ref{fig:336model} (a).
\item Kozma--Szirmai Case: This case represents the extremal case where one horoball is the maximal 
permissible inside the fundamental domain in the sense that it is tangent to the face opposite its center in the tetrahedral cell. 
The remaining three horoballs are smaller but of the same type, but only tangent to the larger horoball on the boundary of the fundamental domain. 
Here $s_0=\frac{1}{2}$ and $s_1=0.142857\dots$. See Fig \ref{fig:336model} (b).
\end{enumerate}

The equations of the horoballs centered at the vertices $E_0$ and $E_1$ in projective coordinates with respect to the two parameters $s_1$ and $s_2$ are given by
$$B_0(\bb_0)=\left( 1, \sqrt{\frac{1- s_0}{2}} \sin\theta \cos \phi, \sqrt{\frac{1-s_0}{2}} \sin \theta  \sin \phi,\frac{1+s_0}{2} + \frac{1-s_0}{2} \cos \theta \right),$$
$$B_1(\bb_1)=\left(1, \sqrt{\frac{1-s_1}{2}} \sin \theta  \cos \phi, \frac{1+s_1}{2} + \frac{1-s_1}{2} \sin\theta \sin\phi ,\sqrt{\frac{1- s_1}{2}} \cos\theta\right).$$
The remaining two horoballs at $E_2$ and $E_3$ are found using rotations 
\begin{equation*}
B_2(\bb_2)=\bb_1 \cdot\left(
\begin{array}{cccc}
1 & 0 & 0 & 0 \\
0 & -\frac{1}{2} & -\frac{\sqrt{3}}{2} & 0 \\
0 & \frac{\sqrt{3}}{2} & -\frac{1}{2} & 0 \\
0 & 0 & 0 & 1 \\
\end{array}
\right),~
B_3(\bb_3)=\bb_1 \cdot\left(
\begin{array}{cccc}
1 & 0 & 0 & 0 \\
0 &  -\frac{1}{2} & \frac{\sqrt{3}}{2} & 0 \\
0 &-\frac{\sqrt{3}}{2} & -\frac{1}{2} & 0 \\
0 & 0 & 0 & 1 \\
\end{array}
\right).
\end{equation*}

\subsubsection{Images in the Projective Model of $\mathbb{H}^3$}

\begin{figure}[p]
\begin{center}
\caption{The two optimal horoball packings of the tiling $\{3,3,6\}$ with packings generated by 3, 4, 5, and 6 reflections (crowns, or layers about the base cube) of horoballs about the fundamental domain.}
\begin{tabular}{|c|c|c|}
\hline
Crowns & Kozma--Szirmai case & B\"or\"oczky--Florian case \\
\hline
3&
\includegraphics[height=50mm]{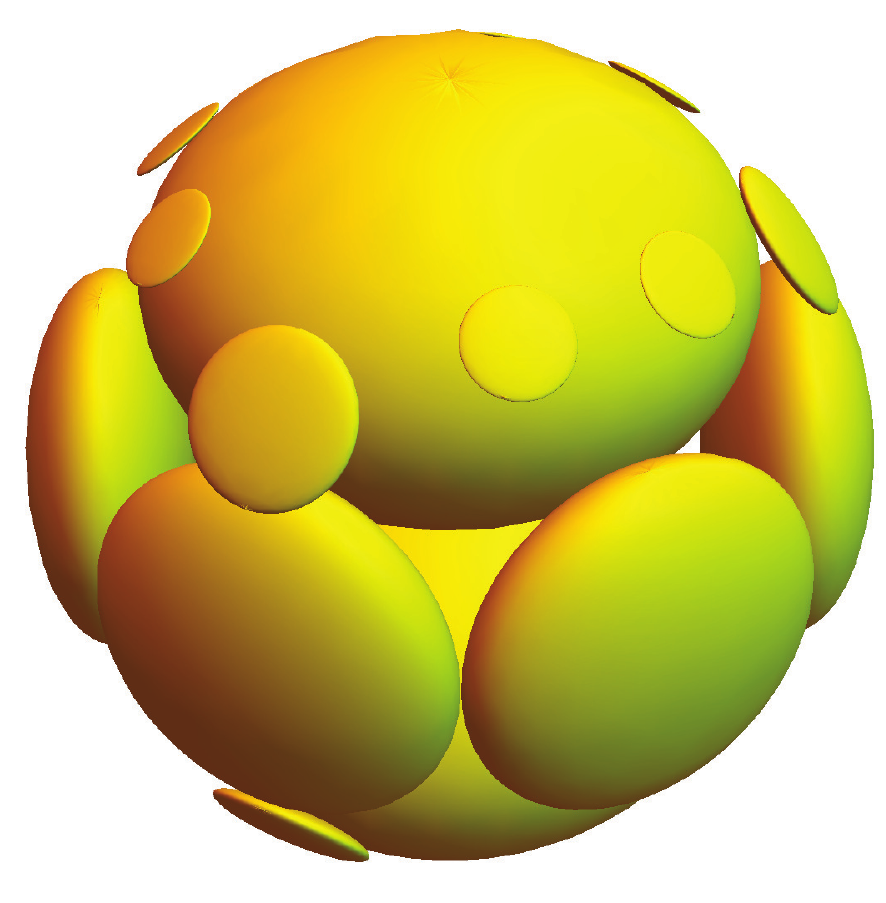}&
\includegraphics[height=50mm]{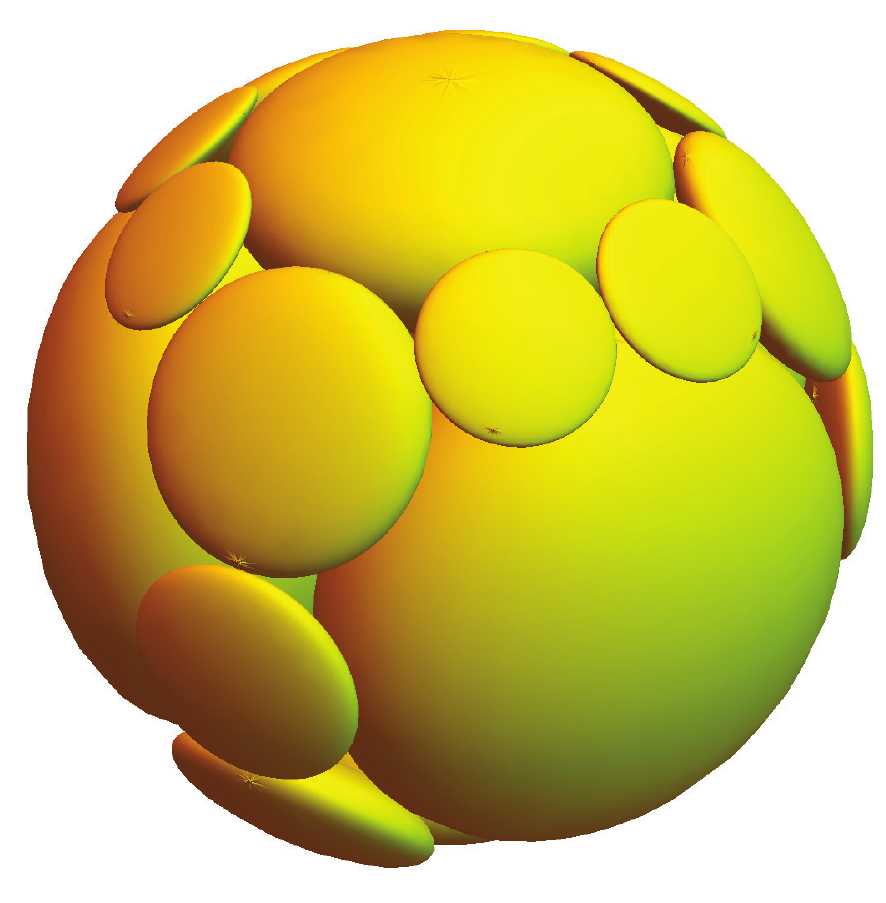}\\
\hline
4&
\includegraphics[height=50mm]{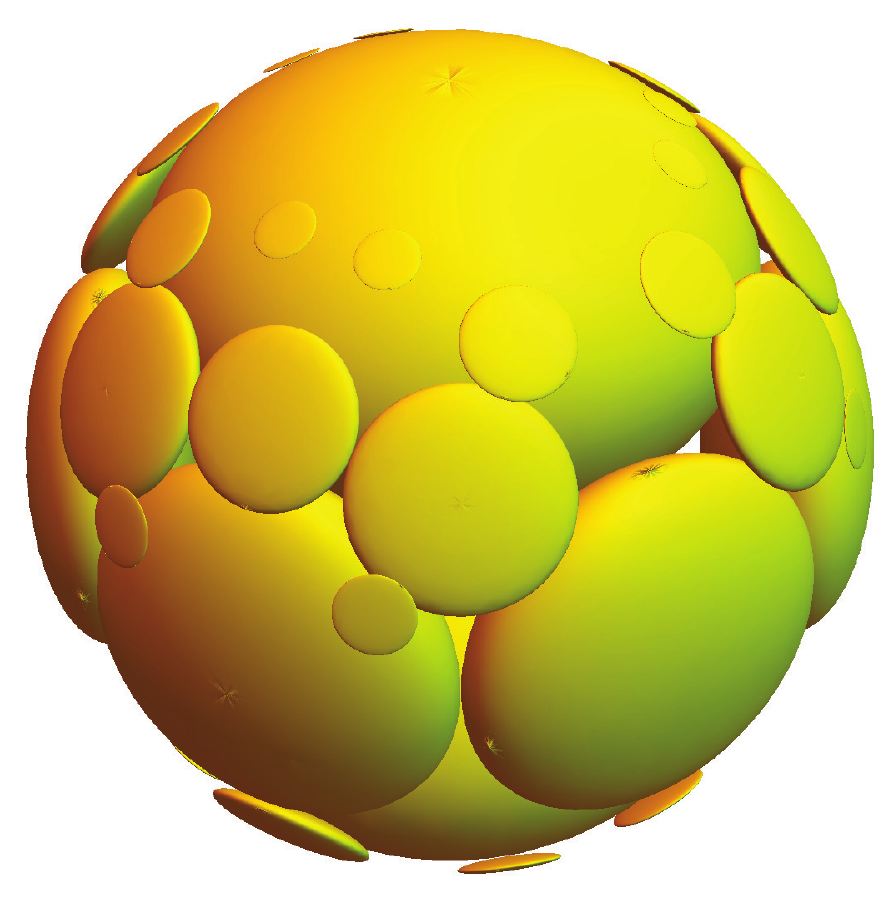}&
\includegraphics[height=50mm]{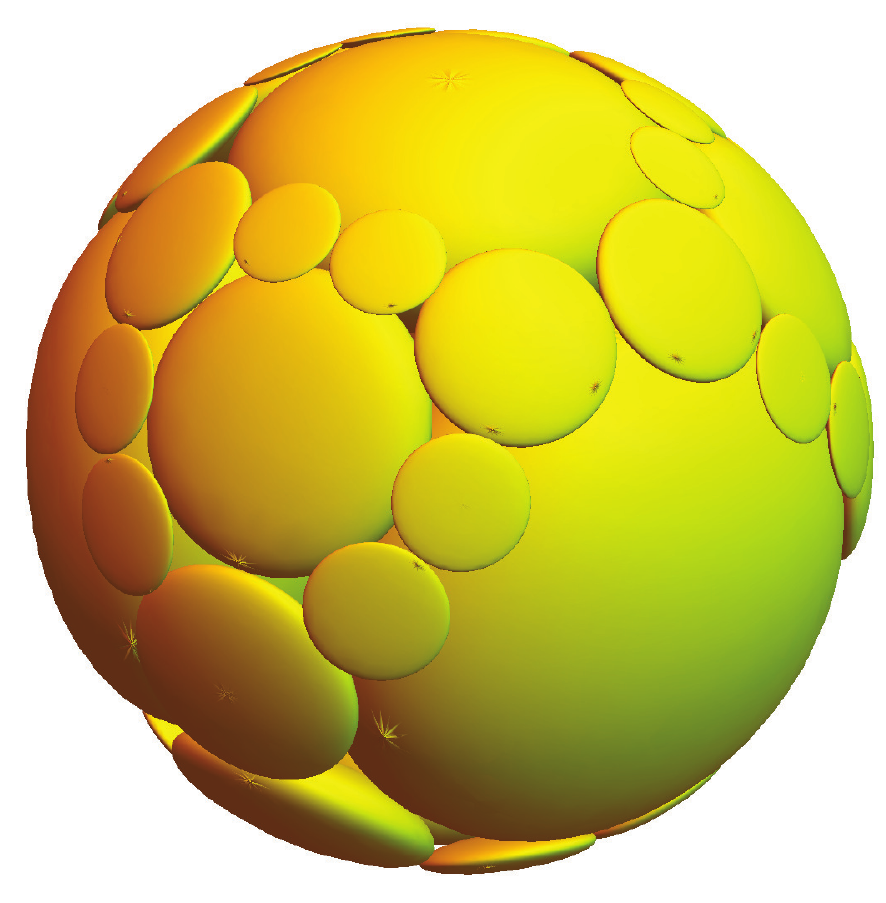}\\
\hline
5&
\includegraphics[height=50mm]{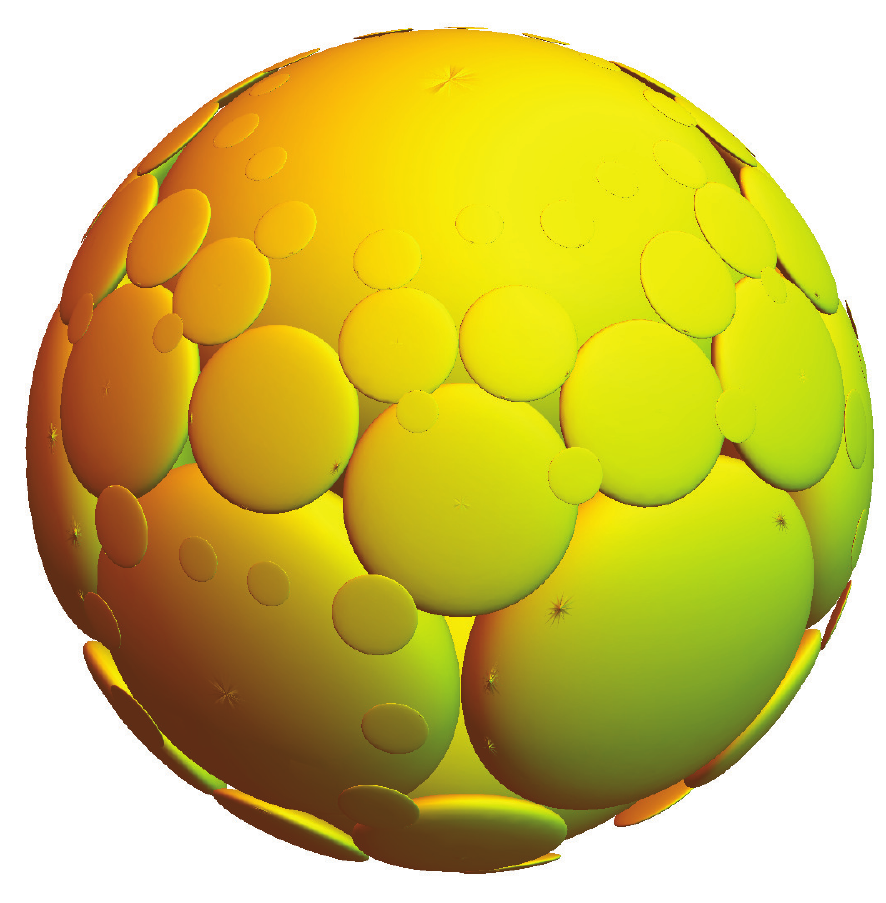}&
\includegraphics[height=50mm]{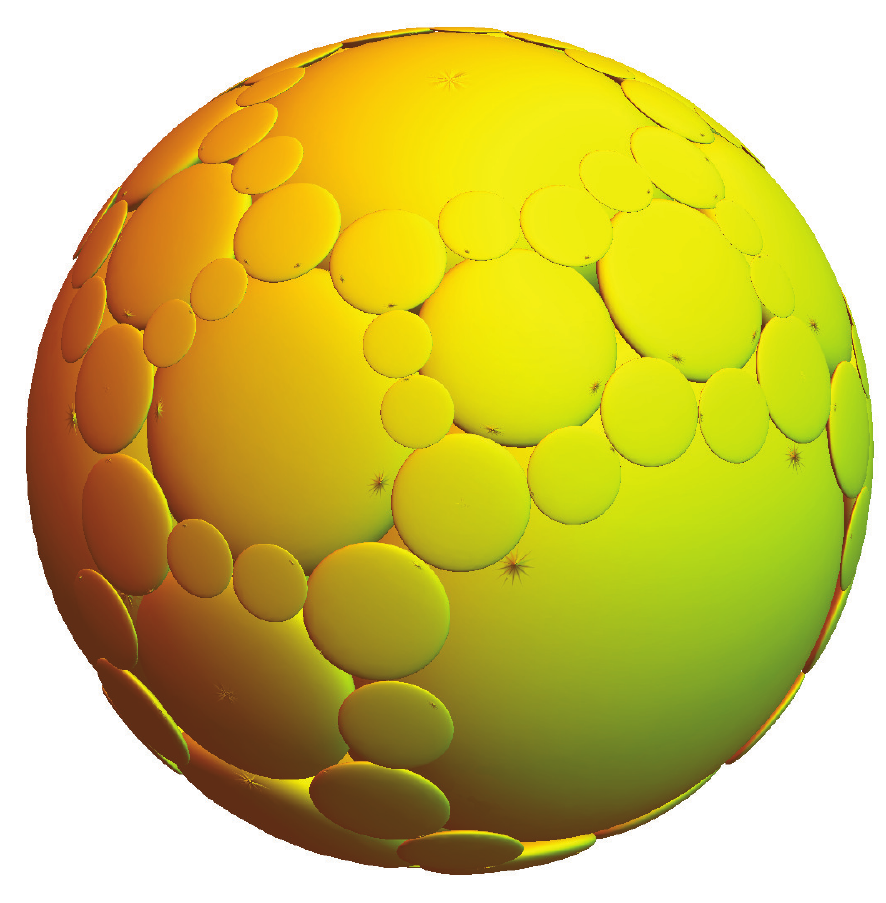}\\
\hline
6&
\includegraphics[height=50mm]{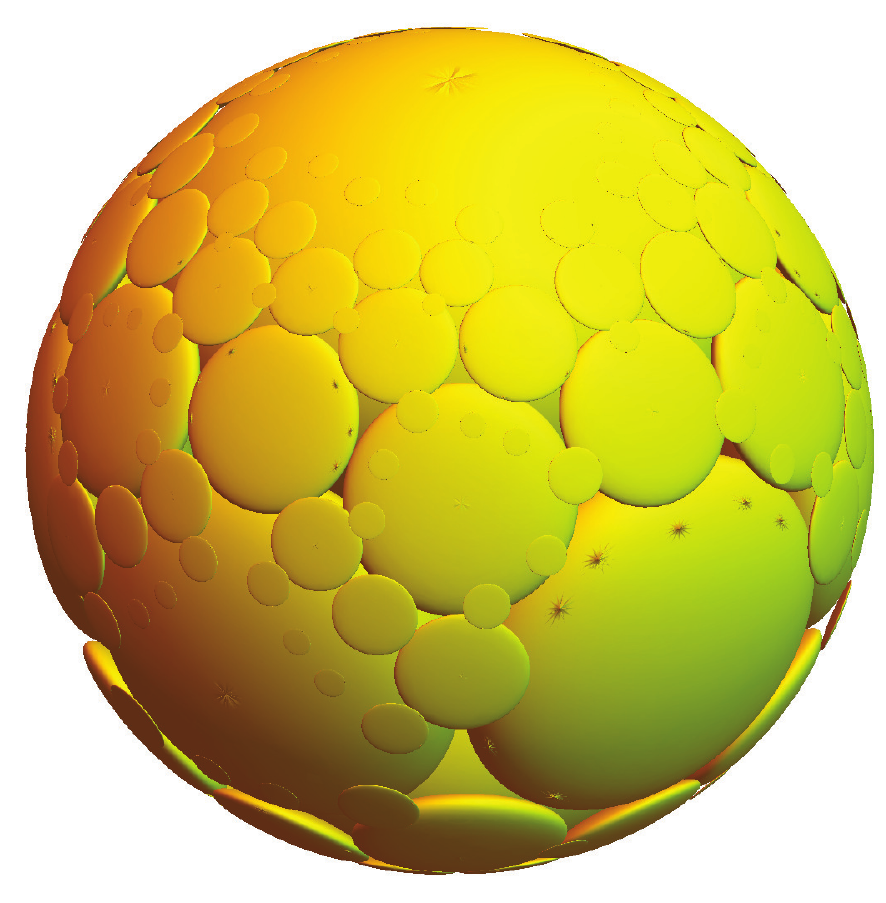}&
\includegraphics[height=50mm]{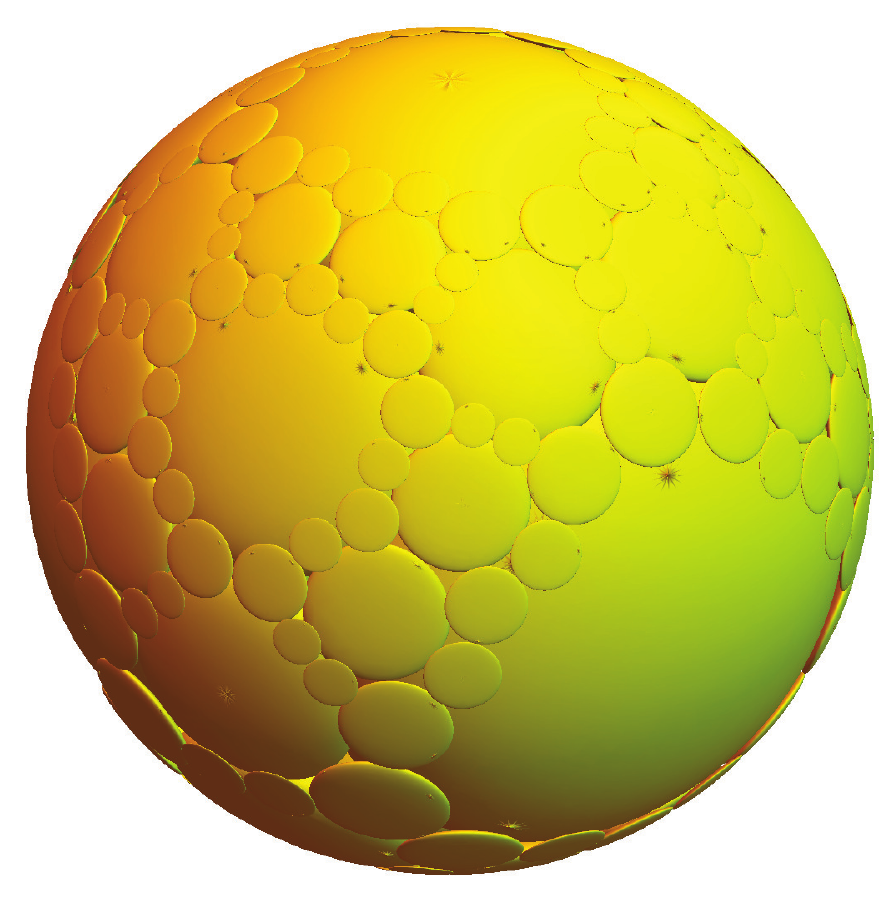}\\
\hline
\end{tabular}
\label{fig:336}
\end{center}
\end{figure}

Applying the generators of the Coxeter group $\bT$ to the fundamental domain we extend the above two packings to $\mathbb{H}^3$, 
and plot the result of the first few iterations in Figure \ref{fig:336}. This is well defined as the regular packings are invariant under the Coxeter group $\bT$.
Figure \ref{fig:336} shows the successive crowns or layers, corresponding to the number of applications of the generators to the base horoballs in the fundamental domain in Figure \ref{fig:336model}.

The figures show the two optimal regular horoball packings of hyperbolic space that arise in the tessellation 
with Schl\"afli symbol \{3,3,6\}, in the Beltrami-Klein model. 
Notice the fractal structure that arises from the embedding into Euclidean 3-space of the packing. 
The balls may appear to have different size with respect to the Euclidean metric on the embedding, but with respect to the hyperbolic metric the balls are congruent. 
The Beltrami-Klein model is not conformal (does not preserve angles) so the balls appear as ellipsoids. 
For proofs of the optimality of these packings see \cite{KSz}. Note that in the $x-y$ slice of the tetrahedral case we have the $\{\infty,3\}$ packing configuration of the hyperbolic plane, see Section 2.

\subsection{$\{4,3,6\}$ Cubic Tiling}
In analogy to the tetrahedral case, we fix the vertex set $\{E_1, E_2, \dots , E_8\}$ of the fundamental domain of the cubic lattice generated by $\bC$ in the Beltrami--Klein model.
\begin{align*}
E_0 &= (1, 0, 0, 1),& 
E_1 &= (1, -\sqrt{2}/\sqrt{3}, \sqrt{2}/3, 1/3), \\
E_2 &= (1, \sqrt{2}/\sqrt{3}, \sqrt{2}/3, 1/3), &
E_3 &= (1, 0, -(2 \sqrt{2})/3, 1/3), \\
E_4 &= (1, 0, (2 \sqrt{2}/3, -1/3),  &
E_5 &= (1, -\sqrt{2}/\sqrt{3}, -\sqrt{2}/3, -1/3),\\ 
E_6 &= (1, \sqrt{2}/\sqrt{3}, -\sqrt{2}/3, -1/3), &
E_7 &= (1, 0, 0, -1). 
\end{align*}
Table \ref{table:data_436} summarizes the data of the fundamental domain of $\mathbf{C}$ used to compute the 
group generators $\bg_i$.
Solving the analogous linear system as above we obtain the generators $\bg_i$ of the group $\mathbf{C}$ and the vertices of the corresponding cubic tiling.  
The group generators are given by
\begin{align*}
\bg_1&
=\left(
\begin{array}{cccc}
 2 & 0 & -\sqrt{2} & -1 \\
 0 & 1 & 0 & 0 \\
 \sqrt{2} & 0 & -1 & -\sqrt{2} \\
 1 & 0 & -\sqrt{2} & 0 \\
\end{array}
\right),
&
\bg_2&=\left(
\begin{array}{cccc}
 -2 & -\sqrt{\frac{3}{2}} & -\frac{1}{\sqrt{2}} & 1 \\
 \sqrt{\frac{3}{2}} & \frac{1}{2} & \frac{\sqrt{3}}{2} & -\sqrt{\frac{3}{2}} \\
 \frac{1}{\sqrt{2}} & \frac{\sqrt{3}}{2} & -\frac{1}{2} & -\frac{1}{\sqrt{2}} \\
 -1 & -\sqrt{\frac{3}{2}} & -\frac{1}{\sqrt{2}} & 0 \\
\end{array}
\right),
\\
\bg_3&=\left(
\begin{array}{cccc}
 2 & -\sqrt{\frac{3}{2}} & \frac{1}{\sqrt{2}} & -1 \\
 \sqrt{\frac{3}{2}} & -\frac{1}{2} & \frac{\sqrt{3}}{2} & -\sqrt{\frac{3}{2}} \\
 -\frac{1}{\sqrt{2}} & \frac{\sqrt{3}}{2} & \frac{1}{2} & \frac{1}{\sqrt{2}} \\
 1 & -\sqrt{\frac{3}{2}} & \frac{1}{\sqrt{2}} & 0 \\
\end{array}
\right),
&
\bg_4&=\left(
\begin{array}{cccc}
 2 & \sqrt{\frac{3}{2}} & -\frac{1}{\sqrt{2}} & 1 \\
 -\sqrt{\frac{3}{2}} & -\frac{1}{2} & \frac{\sqrt{3}}{2} & -\sqrt{\frac{3}{2}} \\
 \frac{1}{\sqrt{2}} & \frac{\sqrt{3}}{2} & \frac{1}{2} & \frac{1}{\sqrt{2}} \\
 -1 & -\sqrt{\frac{3}{2}} & \frac{1}{\sqrt{2}} & 0 \\
\end{array}
\right),
\\
\bg_5&=\left(
\begin{array}{cccc}
 -2 & \sqrt{\frac{3}{2}} & \frac{1}{\sqrt{2}} & -1 \\
 -\sqrt{\frac{3}{2}} & \frac{1}{2} & \frac{\sqrt{3}}{2} & -\sqrt{\frac{3}{2}} \\
 -\frac{1}{\sqrt{2}} & \frac{\sqrt{3}}{2} & -\frac{1}{2} & -\frac{1}{\sqrt{2}} \\
 1 & -\sqrt{\frac{3}{2}} & -\frac{1}{\sqrt{2}} & 0 \\
\end{array}
\right),
&
\bg_6&=\left(
\begin{array}{cccc}
 2 & 0 & \sqrt{2} & 1 \\
 0 & 1 & 0 & 0 \\
 -\sqrt{2} & 0 & -1 & -\sqrt{2} \\
 -1 & 0 & -\sqrt{2} & 0 \\
\end{array}
\right).
\end{align*}

\begin{figure}[p]
\begin{center}
\caption{The two optimal horoball configurations of the cubic tiling $\{4,3,6\}$. The black cube is the fundamental domain, and horoballs are centered at each vertex. The outer sphere represents the boundary of $\partial \mathbb{H}^3$ in the Beltrami--Klein model.}
\begin{tabular}{cc}
\includegraphics[height=60mm]{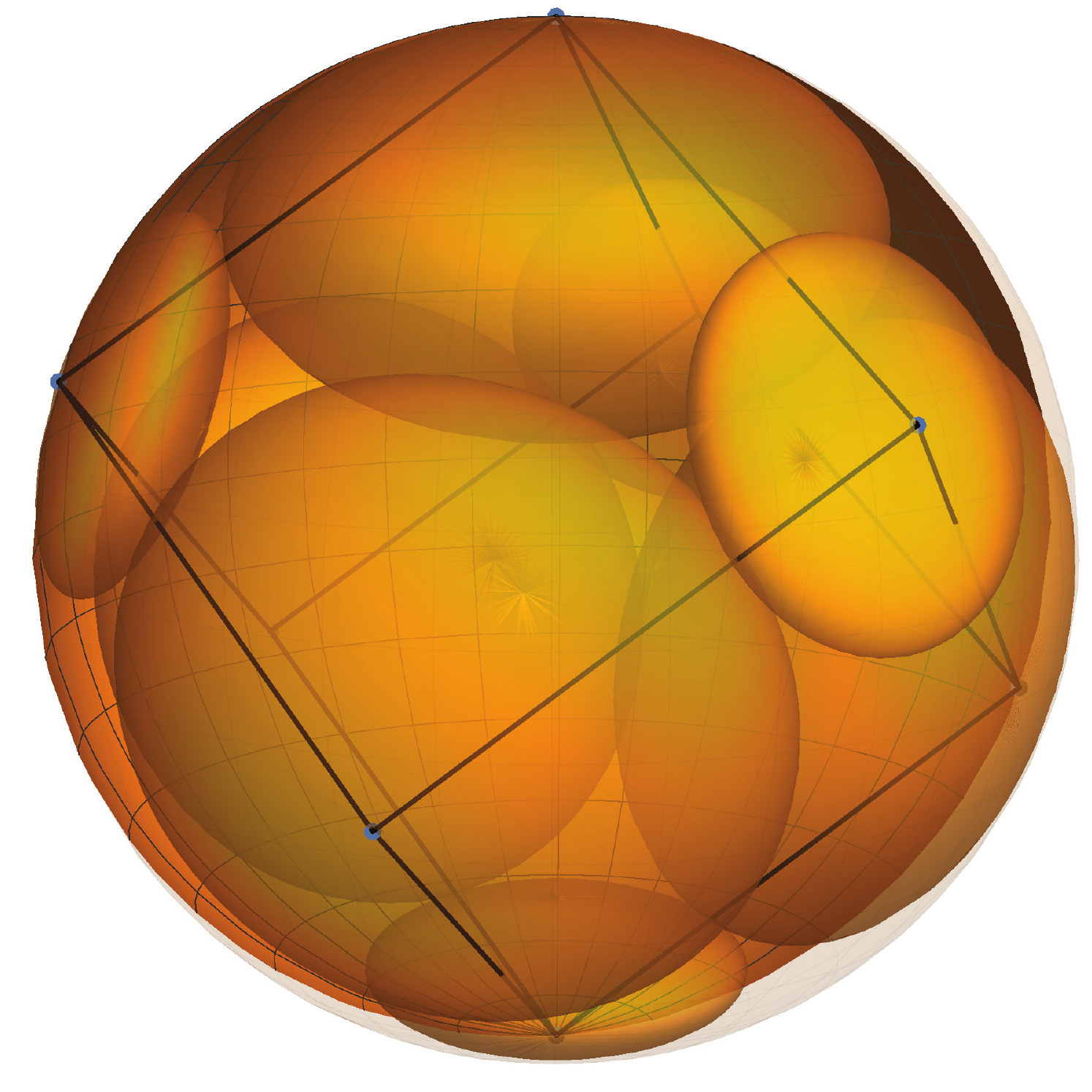} &
\includegraphics[height=60mm]{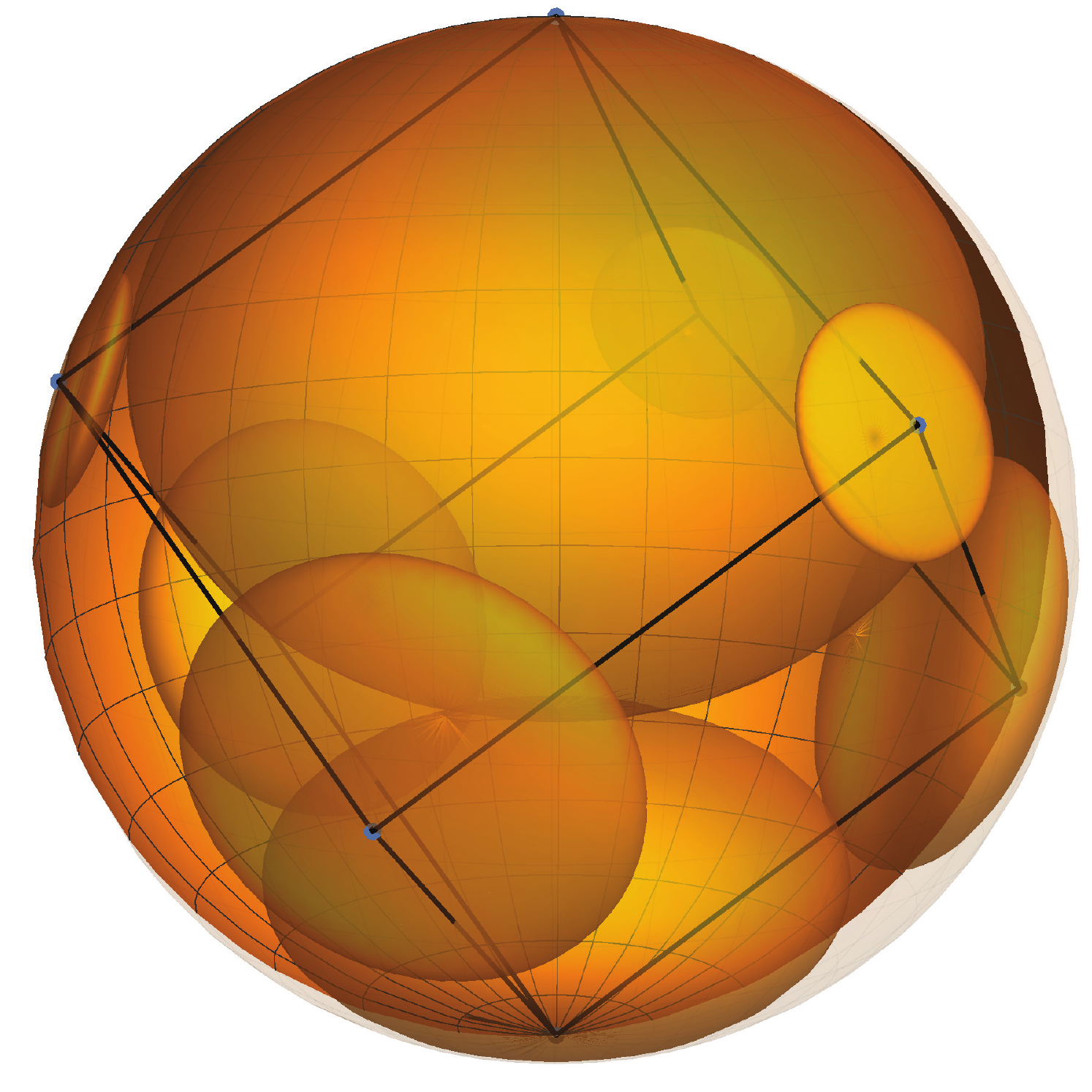} \\
 (a) Maximal case & (b) Balanced case
\end{tabular}
\label{fig:436_model}
\end{center}
\end{figure}

Applying rotations to the polar form horoball equations centered at $(1,0,0,1)$  we find the eight fundamental horospheres with centers at vertices $E_0,\dots E_7$.  
Again, two optimal horoball packings configurations exist that yield the optimal packing density of $0.85328\dots$, see Figure \ref{fig:436_model} (a) and (b). In the first case, there are two horoball types, the larger four are tangent at the midpoints of the facets, in the second, one horoball is the largest admissible in the cubic cell. We refer to the first as the balanced case, and the second as the maximal case. 

We extend these packing from the fundamental domain by symmetries of the Coxeter group $\bC$, see Figure \ref{fig:436} for plots of the two packings with $3$, $4$, and $5$ crowns (layers).

\section{Concluding Remarks}
In this paper we developed a procedure to investigate and visualize the structure of the optimal horoball arrangements 
in $3$-dimensional hyperbolic space.
Figures \ref{fig:336} and \ref{fig:436} show that the optimal packings appear to be structurally distinct based on their contact structures. In general the group isomorphism problem for finitely generated groups is difficult, and goes beyond the scope of this paper. Our figures resemble those of Apollonian gaskets and familiar limit sets of Kleinian groups. 

We note here that similar packings and fractal like images can be derived from the higher dimensional optimal horoball packings described in \cite{KSz14}. We also point out that the existence of multiple optima or equilibrium states for packings in $\mathbb{H}^3$ may have nontrivial consequences for the physical sciences. 
 
 \section{Acknowlegements}
 We would like to thank K\'aroly B\"or\"oczky for his valuable conversations during the L\'aszl\'o Fejes T\'oth Centennial conference in Budapest, which lead to the figures presented in this paper.




\begin{table}[p]
	\begin{tabular}{|l|l|l|l|}
		 \hline
		 Facet & Normal & Footpoint & Reflection\\ 
		 \hline
		 $E_1E_2E_3$ & $\Bu_0=(1,0,0,1)$ & $\textbf{Pr}_{\Bu_0}(E_0)=(1,0,0,0)$ & $\bg_1(E_0)=(1,0,0,-1)$\\
		 $E_0E_2E_3$ & $\Bu_1=(1,0,2,-1)$ & $\textbf{Pr}_{\Bu_1}(E_1)=(1,0,-\frac{2}{7},\frac{3}{7})$ & $\bg_2(E_1)=(1,0,-\frac{4}{5},\frac{3}{5})$\\
		 $E_0E_1E_3$ & $\Bu_2=(-1 , -\sqrt{3} , 1, 1)$ & $\textbf{Pr}_{\Bu_2}(E_2)=(1,-\frac{\sqrt{3}}{7},\frac{1}{7},\frac{3}{7})$ & $\bg_3(E_2)=(1,-\frac{2 \sqrt{3}}{5},\frac{2}{5},\frac{3}{5})$\\
		 $E_0E_1E_2$ & $\Bu_3=(1 , -\sqrt{3} , -1, -1)$ & $\textbf{Pr}_{\Bu_3}(E_3)=(1,\frac{\sqrt{3}}{7},\frac{1}{7},\frac{3}{7})$ & $\bg_4(E_3)=(1,\frac{2 \sqrt{3}}{5},\frac{2}{5},\frac{3}{5})$\\
		 \hline 
	\end{tabular}
	\caption{Metric data used to find the generators of group $\bC$.}
	\label{table:data_336}
\end{table}


\begin{table}[p]
	\begin{tabular}{|l|l|l|l|}
		 \hline
		 Facet & Normal & Footpoint & Reflection\\ 
		 \hline
		 $E_0E_1E_2E_4$ & $\Bu_1=(1,0,-\sqrt{2},-1)$ & $\textbf{Pr}_{\Bu_1}(E_7)=(1,0,\frac{1}{\sqrt{2}},0)$ & $\bg_1.E_7=(1,0,\frac{2 \sqrt{2}}{3},\frac{1}{3})$\\
		 $E_0E_1E_3E_5$ & $\Bu_2=(-2,-\sqrt{6},-\sqrt{2},2)$ & $\textbf{Pr}_{\Bu_2}(E_7)=(1,-\frac{\sqrt{3}}{2\sqrt{2}},-\frac{1}{2 \sqrt{2}},0)$ & $\bg_2(E_7)=(1,-\sqrt{\frac{2}{3}},-\frac{\sqrt{2}}{3},\frac{1}{3})$\\
		 $E_0E_2E_3E_6$ & $\Bu_3=(2,-\sqrt{6},\sqrt{2},-2)$ & $\textbf{Pr}_{\Bu_3}(E_7)=(1,\frac{\sqrt{3}}{2\sqrt{2}},-\frac{1}{2 \sqrt{2}},0)$ & $\bg_3(E_7)=(1,\sqrt{\frac{2}{3}},-\frac{\sqrt{2}}{3},\frac{1}{3})$\\
		 $E_7E_1E_4E_5$ & $\Bu_4=(2,\sqrt{6},-\sqrt{2},2)$ & $\textbf{Pr}_{\Bu_4}(E_0)=(1,-\frac{\sqrt{3}}{2\sqrt{2}},\frac{1}{2 \sqrt{2}},0)$ & $\bg_4(E_0)=(1,-\sqrt{\frac{2}{3}},\frac{\sqrt{2}}{3},-\frac{1}{3})$\\
		 $E_7E_2E_4E_6$ & $\Bu_5=(-2,\sqrt{6},\sqrt{2},-2)$ & $\textbf{Pr}_{\Bu_5}(E_0)=(1,\frac{\sqrt{3}}{2\sqrt{2}},\frac{1}{2 \sqrt{2}},0)$ & $\bg_5(E_0)=(1-\sqrt{\frac{2}{3}},\frac{\sqrt{2}}{3},-\frac{1}{3})$\\
		 $E_7E_3E_5E_6$ & $\Bu_6=(1,0,\sqrt{2},1)$ & $\textbf{Pr}_{\Bu_6}(E_0)=(1,0,-\frac{1}{\sqrt{2}},0)$ & $\bg_6(E_0)=(1,0,-\frac{2 \sqrt{2}}{3},-\frac{1}{3})$\\
		 \hline 
	\end{tabular}
	\caption{Metric data used to find the generators of group $\bT$.}
	\label{table:data_436}
\end{table}

\end{document}